\documentclass[12pt]{article}
\usepackage{amsmath}
\usepackage{amssymb}
\usepackage{graphicx}
\usepackage{amsfonts}
\usepackage{amsopn}
\usepackage{amsthm}

\def\Ca{C^{1,\a}_a(\Sn)}

\def\Reals{{\mathbb R}}

\def\Rn{{\mathbb R}^{n+1}}
\def\R3{{\mathbb R}^3}
\def\Rnn{{\Reals}^{n}}

\def\Rb{\bar{R}}
\def\xb{\bar{x}}

\def\t{{\bf t}}

\def\calA{{\mathcal A}}

\def\calK{{\mathcal K}}

\def\calO{{\mathcal O}}

\def\Sn{{\mathbb S}^n}
\def\R{{\mathcal R}}
\def\M{{\mathcal M}^n}

\def\aa{{\bf a}}
\newcommand{\ep}{\epsilon}
\renewcommand{\O}{\Omega}
\newcommand{\p}{\rho}
\renewcommand{\t}{\tau}

\newcommand{\n}{\nabla}

\def\ka{{\kappa}}

\newcommand{\ox}{\bar {x}}
\newcommand{\ov}{\bar {v}}

\newcommand{\ow}{\bar {w}}

\newcommand{\eh}{\hat{e}}

\newcommand{\gba}{\bar {g}}

\renewcommand{\o}{\omega}
\newcommand{\s}{\sigma}

\newcommand{\gR}{{\gamma_R}}

\renewcommand{\a}{\alpha}

\newcommand{\D}{\Delta}
\newcommand{\Dh}{\hat{\Delta}}

\newcommand{\del}{\partial}

\newcommand{\ra}{\rangle}
\newcommand{\G}{\Gamma}
\newcommand{\lm}{\lambda}

\def\calA{{\mathcal A}}

\def\calK{{\mathcal K}}

\newcommand{\Ub}{\bar{U}}
\renewcommand{\O}{\Omega}
\newcommand{\la}{\langle}

\newtheorem{theorem}{Theorem}

\newtheorem{lemma}[theorem]{Lemma}
\newtheorem{proposition}[theorem]{Proposition}
\newtheorem{remark}[theorem]{Remark}
\newtheorem*{remark*}{Remark}
\title{Kummer configurations and $S_m-$reflector problems: Hypersurfaces in $\Rn$ with given mean intensity
}
\author{Vladimir I. Oliker\thanks{The research of the author
was partially supported by National Science Foundation grant DMS-04-05622.
Part of research for this paper
was performed while the author was a guest of the Technical University of
Berlin during December 2002 - January 2003. The final version of this paper was written during another visit to the TU-Berlin in June-July of 2008. The author is grateful to
the Fachbereich Mathematik of the TU-Berlin, and particularly, to Professor
Udo Simon for the warm hospitality.}\\
Department of Mathematics and Computer Science,\\
 Emory University, Atlanta, Georgia\\
oliker@mathcs.emory.edu}
\date{}
\begin{document}
\maketitle
\pagenumbering{arabic}
\setcounter{section}{0}
\setcounter{subsection}{0}
\begin{abstract} For a congruence of straight lines defined by a hypersurface in $\Rn,~n \geq 1,$ and a field
of reflected directions created by a point source we define the notion of intensity in a tangent direction and
introduce elementary symmetric functions $S_m,~m=1, 2,...,n,$ of {\it principal intensities}.
The problem of existence and uniqueness of a closed hypersurface with prescribed $S_n$ is the
``reflector problem'' extensively studied in recent years. In this paper we formulate and give sufficient
conditions for solvability of an analogous problem in which the mean intensity $S_1$ is a given
function.
\end{abstract}

\section{Introduction}\label{intro}
Fix a Cartesian coordinate system in $\Rn, ~n \geq1,$ with
the origin $\calO$
and
let $\Sn$ be the unit sphere centered at $\calO$. It was shown in \cite{galo:04} that among measure-preserving maps of $\Sn$ onto  itself transferring two given positive Borel measures  into each other there exists a uniquely defined map which is optimal against
the cost function with density $-\log(1-\la x, y\ra),~x,y \in \Sn$; here $\la \cdot,\cdot \ra$ denotes the inner product in $\Rn$. This result was established under
 quite mild assumptions on the given measures. Even prior to \cite{galo:04} it was already shown in \cite{glimm-oliker-refl:03} (see Theorem 3.4 and the note following its proof), that such optimal map is generated by a closed
convex hypersurface $R$ in $\Rn$  which is star-shaped relative to $\calO$ and acts as a reflector for light rays emanating from $\calO$. At smooth points of $R$ this optimal map is given by
\begin{equation}\label{snell}
y=\gR(x) = x - 2 \la x,N(x)\ra N(x),
\end{equation}
where $N$ denotes the unit normal field on $R$. This, of course, is the classical law of reflection. (At nonsmooth points of $R$ there is an appropriate generalization of
the map $\gR$.)

Since $R$ is convex, it is almost everywhere of class $C^2$ and one may consider
the Jacobian $J(\gR)$ defined almost everywhere on $\Sn$. Then, from the geometric point of view, the $|\det J(\gR)|$ is the quotient
of densities of the volume forms defining the two given Borel measures. This relation
leads to a (possibly, degenerate) second order elliptic partial differential equation (PDE) of Monge-Amp\`{e}re type on $\Sn$ for the radial  function defining $R$ \cite{OW}. (The radial function is defined at the beginning of section \ref{radial} below.) Thus, the geometric problem of finding a hypersurface $R$ such that the map (\ref{snell}) transfers two given volume forms into each other requires solution of the corresponding fully nonlinear PDE. This problem is usually referred to as the {\it reflector problem}.  Existence of a weak solution to the reflector  problem was shown in \cite{CO} for surfaces in
${\mathbb R}^3$ but the proof is valid verbatim for hypersurfaces in $\Rn$.  Uniqueness was shown
in \cite{Pegfei_guan_wang:jdg98}, \cite{glimm-oliker-refl:03} and  regularity was studied in \cite{Pegfei_guan_wang:jdg98}, \cite{CGH:04}, \cite{loeper:06}.

The described results suggest that the map $\gR$ is interesting from several points of view and deserves further investigation. Indeed, in this paper we show that the reflector problem  is only one of a series of semilinear and fully nonlinear geometric problems connected with the map $\gR$.

A very natural geometrical framework for studying the map $\gR$ is the Kummer configuration
considered by E. Kummer in 1860  in his paper \cite{kummer} on  congruences of straight lines in ${\mathbb R}^3$.
Congruences of straight lines arise naturally in geometrical optics and optimal mass transport in $\Rnn$ and were considered (in ${\mathbb R}^3$) already in the 18-th century by G. Monge and in the early part of 19-th century by E. Malus and W.R. Hamilton. In \cite{kummer} Kummer    defines a congruence of straight lines in ${\mathbb R}^3$  by points on a given surface (base) and a set of direction vectors. This pair, the surface and the vector field, is referred to as a ``{\it Kummer configuration}''; see, \cite{kagan}, v. 2, ch. 17. For such  a congruence Kummer introduced notions analogous to the first and second fundamental forms (the latter is not necessarily symmetric!) and studied  its properties which can be described  using these forms.
In  the years subsequent to the publication of \cite{kummer}, the dependence of the second fundamental form  of Kummer on the base surface was considered by geometers as a deficiency and theories avoiding such dependence were developed \cite{Finikov:1950}.

In this paper we treat the hypersurface $R$  as a reflector and the reflected rays defined by the map $\gR$ as a congruence of straight lines, that is, we have a special case of
a Kummer configuration $(R,\gR)$. This point of view is our starting point, even though the definitions and objectives here are different from those of Kummer.

The paper is organized as follows. In section \ref{radial} we describe the class of hypersurfaces in $\Rn$
for which the map $\gR$ is  studied  and derive various local formulas. In section \ref{kummer}
we introduce the notion of {\it intensity in direction of a curve} and show that in principal directions the {\it principal} intensities
are the real eigenvalues of a certain quadratic differential form analogous to the second fundamental form in classical differential geometry.  In the same section we introduce the elementary symmetric functions $S_m$ of principal intensities; here $m$ is an integer, $1 \leq m \leq n$. The problem of finding the optimal map
described in the first paragraph of this introduction, that is, the reflector problem, corresponds to $m=n$.
In section \ref{ME_1} we establish existence and uniqueness of solutions to the $S_1$-reflector problem.
We intend to present solutions to other reflector problems in a separate publication.

The author is indebted to the referee for reading carefully the manuscript and for useful comments.

\section{\protect\bf Reflectors defined by radial functions}\label{radial}
In this section our considerations are local. Let $x=x(u)\equiv x(u^1,...,u^n)$ be a smooth local parametrization of $\Sn$. Let $R$ be a hypersurface in $\Rn$ which is a graph
over some domain $\o \subset S^n$ of a function $\p: \o \rightarrow (0,\infty), ~\p \in C^2(\o)$. Such $R$ can be defined by the position vector
$r(x)= \p (x)x,~ x \in \o$. (In this paper  $x\in \Sn$ is treated as a point in $\Sn$ and also as a unit vector in $\Rn$.)
The function $\p$ is called
the {\it radial function of $R$}.  Obviously, the map $r:\o \rightarrow \Rn$ is an embedding. The set of all such hypersurfaces in $\Rn$ is denoted
by $\M$.
If we need to indicate the domain $\o$ we write $\M(\o)$; in particular, if $\o \equiv \Sn$ we write $\M(\Sn)$.
We will study reflecting properties of hypersurfaces in $\M$ and
for brevity refer to them as {\it reflectors}.

Denote by $e=e_{ij}du^idu^j$ the standard metric on $\Sn$ induced from $\Rn$. Here and for the rest of the paper the  Latin indices $i,j,k,...$ run over
 the range $1,2, ...,n$ and the summation convention over repeated lower
 and upper indices is in effect. The following notation will be used:
 \[\del_i = \frac{\del}{\del u_i}, ~x_i = \del_ix, ~r_i = \del_i r,~\mbox{etc.}~  \del_{ij} = \frac{\del^2}{\del u_i\del u_j},~x_{ij}=\del_{ij}x, ~r_{ij} = \del_{ij} r,~\mbox{etc}. \]
 The covariant differentiation
 in the metric $e$ is denoted by $\n_i:=\n_{\del_i}$ and similarly $\n_{ij}$,
etc. On functions, $\n_i = \del_i$ and
$\n_{ij} = \del_{ij} - \G_{ij}^k\del_k$, where $\G_{ij}^k$ are the Christoffel symbols of the metric $e$. Put $\n = e^{ij}x_j\del_i$, where $[e^{ij}]=[e_{ij}]^{-1}$, and
$W_{\p}= \sqrt{\p^2 + |\n \p|^2}$.

Let $R\in \M$.  We recall first the expressions for the classical first and second fundamental forms
of $R$ in terms of its radial function \cite{O:curv}.  The coefficients $g_{ij}$ of the first fundamental form $g$ of $R$, the elements of the
inverse matrix $[g^{ij}] = [g_{ij}]^{-1}$, and the determinant of $[g_{ij}]$ are, respectively,
\begin{eqnarray} \label{form1}
 g_{ij} = \la  r_i,  r_j \ra =
 \p_i \p_j + \p^2 e_{ij},~~ g^{ij}
= \frac{1}{\p^2} \left (e^{ij} - \frac{\p^i\p^j}{W_{\p}^2} \right ),\nonumber\\~\det[g_{ij}] = \p^{2n-2}W_{\p}^2\det[e_{ij}],
\end{eqnarray}
where $\p^i = e^{ik}\p_k.$
The unit normal field $N$ on $R$ is given by
\begin{equation} \label{normal}
N = \frac{\p x - \n \p}{W_{\p}}.
\end{equation}
The coefficients of the second fundamental form of $R$ are given by
\begin{equation} \label{formb}
b_{ij} = -\la r_j, N_i \ra =
\frac{\p \n_{ij}\p -2\p_i\p_j - \p^2e_{ij}}{W_{\p}}.
\end{equation}
Since $r_i = \p_ix+\p x_i$, it follows from (\ref{snell}) that $\la r_i, \gR \ra = \p_i$. Differentiating, we obtain
\[ \la r_{ij}, \gR \ra + \la r_i, {\gR}_j \ra= \p_{ij}.\]
This implies that
\begin{equation}\label{symm1}
\la r_i, {\gR}_j \ra = \la r_j, {\gR}_i \ra.
\end{equation}

We will need explicit expressions of
$\la r_i, \gR_j \ra$ and $\la \gR_i, \gR_j \ra$ in terms of $\p$ and its derivatives.
To determine
$\la r_i, \gR_j \ra$,
differentiate (\ref{snell}) and take the inner product of the
result with $r_i$. Then
\[ \la r_i,\gR_j \ra
= \p e_{ij}   +2 \la x, N \ra b_{ij}.\]
Put
\begin{equation}\label{ka}
\ka_{ij}:=-\frac{\la r_i, \gR_j \ra}{\p}~~~\mbox{and}~~~\eh_{ij}:=\la \gR_i, \gR_j \ra.
\end{equation}
Noting that
$\la x,N \ra = \p/W_{\p}$, we
get
\begin{equation} \label{H}
 -\ka_{ij}= e_{ij} + \frac{2}{W_{\p}}b_{ij}.
\end{equation}
For reasons which will become clear in a moment the quadratic differential form $\ka=\ka_{ij}du^idu^j$ will be referred to as the {\bf intensity} form of the congruence $(R,\gR)$. Its geometric meaning will also be described below.

Next, we derive an expression for $\eh_{ij}$ in terms of $\p$. Note first that because for each $x\in \o$ the vectors $r_1(x),...,r_n(x), N(x)$  form a basis of $R^{n+1}$ we have
\begin{equation} \label{gidec}
\gR_i = \la \gR_i, r_s \ra g^{sk}r_k + \la \gR_i, N \ra N.
\end{equation}
Using (\ref{snell}), (\ref{normal}), the equations of Weingarten $N_i = -b_{ij}g^{jk}r_k$ and noting that by (\ref{form1})
$\p_k g^{kj} = \frac{\p^j}{W_{\p}^2}$, we get
\begin{equation} \label{giu1}
\la \gR_i, N \ra= -\la x_i,N \ra - 2 \la x,N_i \ra= \frac{\p_i}{W_{\p}} +
2b_{ij}g^{jk}\p_k=-\frac{\p^j}{W_{\p}}\ka_{ji}.
\end{equation}
It follows from (\ref{H}), (\ref{gidec}), (\ref{giu1}) and (\ref{form1}) that
\begin{equation} \label{gigj}
\eh_{ij}=\la \gR_i, \gR_j \ra  =  \ka_{ik}
\left (e^{kl} -\frac{\p^k\p^l}{W_{\p}^2}\right )
 \ka_{lj} + \ka_{ik} \frac{\p^k\p^l}{W_{\p}^2}
\ka_{lj}= \ka_{ik}e^{kl}\ka_{lj}.
\end{equation}
\section{The Kummer configuration, the intensity form and the $S_m$-reflector problem} \label{kummer}
It is clear from the discussion in the Introduction that the pair $(R,\gR)$ forms a Kummer configuration
 with $R$ as the base hypersurface and $\gR$ defining the directions of reflected rays.
In geometrical optics
the quantity
\begin{equation} \label{density}
|J(\gR(x))|=\frac{\sqrt{\det [\la \gR_i(x),\gR_j(x) \ra]}}{\sqrt{\det [e_{ij}(x)]}}
\end{equation}
is called the {\bf intensity} (or, more accurately, {\bf the relative intensity}) in the reflected direction $\gR(x)$ \cite{W}.
This is a very important quantity
characterizing reflecting properties of  the hypersurface $R$. Assume that the density of the distribution of the light rays emanating
from $\calO$ is given by some function $g(x), ~x \in \o \subset \Sn$. Then the role of  $R$ is to redistribute the energy from the source $\calO$ so that the reflected rays have directions defined by some given region $\O \subset \Sn$ and a prescribed density $f(y),~y\in \O$ \cite{W}. The reflector problem as stated,
for example, in \cite{W},  is to determine such $R$; see \cite{OW} and \cite{Oliker:TiNA} for more details concerning mathematical formulations of this and some related problems.

We clarify now the geometric meaning of the intensity form $\ka$.
Let $\bar{x} \in \o$ and $x(t),~|t-t_0|< \ep$ for some $\ep > 0$, a smooth curve in $\o$ such that $r(x(t_0))=r(\bar{x})$. Denote by $\dot{x}(t)$ the tangent vector to $x(t)$. Define the {\bf intensity
in direction of $x(t)$ at $t=t_0$} as the quotient
\begin{equation}\label{int}
\frac{\sqrt{\eh(\dot{x}(t_0))}}{\sqrt{e(\dot{x}(t_0))}}.
\end{equation}

It follows from (\ref{ka}), (\ref{H}) and (\ref{gigj}) that
\begin{equation}\label{int1}
\mbox{sign}(\ka)\frac{\sqrt{\eh(\dot{x}(t_0))}}{\sqrt{e(\dot{x}(t_0))}}= \frac{\ka (\dot{x}(t_0))}{e(\dot{x}(t_0))},
\end{equation}
where $\ka (\dot{x}(t_0))$ is the value of the form $\ka$ on the tangent vector to $\dot{x}(t_0)$.
Thus, the intensity in direction of $x(t)$ at $t=t_0$ is the rate of change of the angle between
$\gR(x(t))$ and $\gR(x(t_0))$ relative to the change of the angle between
$x(t)$ and $x(t_0)$.

Note that with our choice of the sign,  $\ka$ is positive definite on a unit sphere. This  follows from (\ref{H}) and (\ref{formb}). The coefficients of $\ka$ are clearly invariant with respect to rescaling
$\p \rightarrow \lambda \p$ with $\lambda > 0$. This is
consistent with  the invariance of the map $\gR$ with respect
to homotheties of $R$ with respect to the origin $\calO$. It follows from (\ref{gigj}) that the three forms, $e,\ka$ and $\eh$ are not
independent.

The form $\ka$ has also another  geometric interpretation. Let $x(t)$ be a smooth curve as before. Suppose also that $\gR$ is
a diffeomorphism. Consider the sequence of straight
lines $l(t)$ of directions $\gR(t)$ passing through $r(t)$. Since these
lines are not parallel, for each $t, ~0<|t-t_0| < \ep$, there exists a unique  pair of points $X_0(t) \in l(t_0)$ and $X_t \in l(t)$ realizing
the distance in $\Rn$ between these lines. Denote by $h(t)$ the signed distance
from  $r(t_0)$ to $X_0(t)$ with the ``+'' sign taken if $\la X_0(t)-r(t_0), \gR(t_0)\ra>$ and ``-'' sign otherwise. To calculate the $\lim h(t)$ as $t\rightarrow t_0$,
observe that the segment realizing the distance between  $l(t)$ and $l(t_0)$ is orthogonal to both $\gR(t)$ and $\gR(t_0)$, that is,
$$\la \frac{dr}{dt} + \frac {dh}{dt}\gR+ h\frac{d\gR}{dt},\gR\ra \left. \right | _{t=t_0} =0$$ and
$$\la \frac{dr}{dt} + \frac {dh}{dt}\gR+ h\frac{d\gR}{dt},
\frac{d\gR}{dt}\ra \left. \right | _{t=t_0} =0.$$
Taking into account that $\la \frac{dr}{dt},\gR(t) \ra = \frac{d\p}{dt}$ along the curve $x(t)$ (this follows from (\ref{snell})), the first of this equalities implies $\la\frac{d\p}{dt},\gR(t) \ra = -\frac{dh}{dt}$ at $t=t_0$. The second equality implies (see (\ref{H})) that
\[h(\dot{x}(t_0))= \frac{\p(x(t_0))\ka(\dot{x}(t_0))}{\eh(\dot{x}(t_0))}.\]
A point on the line of direction $\gR(t_0)$ through $r(t_0)$ at the distance $h(\dot{x}(t_0))$ from $r(t_0)$ is called the point of striction. The above formula shows that on each ray there exist a segment filled up with points of striction corresponding to each tangent
direction at $x(t_0)$. Of course, this segments may degenerate into a point.
Note that $\ka(\dot{x}(t_0))$ has a positive sign if the linear segment from the point of reflection to the striction point has the same direction as $\gR(x(t_0))$; otherwise it is negative.

Throughout this paper we adapt the convention that when $d\gR(\dot{x}(t_0))=0$ then
$h(\dot{x}(t_0))=\infty$. In this case, $\ka(\dot{x}(t_0))=0$.


An explicit expression for the $(2,0)$ tensor $\ka$ is obtained using (\ref{H}) and (\ref{formb}):
\begin{equation}\label{kappa0}
\ka = \frac{-\p\n^2\p + 2 \n \p \otimes \n \p +\left[(\p^2-|\n \p|^2)/2\right ]e}{W_{\p}^2/2}~ (\mbox{in}~ \o),
\end{equation}
where $\n^2\p = [\n_{ij} \p],~\n \p \otimes \n \p = [\p_i\p_j]$.

\begin{remark}\label{conformal}
If in (\ref{kappa0}) we make a change $\p = e^{-w}$ then we obtain
\begin{equation}\label{conformal1}
\ka(w) = \frac{\n^2 w +  \n w \otimes \n w +\left[(1-|\n w|^2)/2\right ]e}{(1 + |\n w|^2)/2}.
\end{equation}
The tensor $\kappa(1 + |\n w|^2)/2$ is the Schouten tensor of the metric $e^{-2w}e$ on $\Sn$.
This observation was pointed out to me by M. Gursky \cite{gursky:08}.
\end{remark}

Now we present several examples. For a sphere of radius $R$, using (\ref{formb}) and  (\ref{H}), we obtain $\ka=e$. Similarly, for a piece of a hyperplane in $\Rn$ we have $\ka =-e$.

Consider now an ellipsoid of revolution $E$  with one focus at $\calO$ and axis of revolution passing through both foci. Denote by $\aa$ the second focus.
Using the expressions for $H_{ij}$ in section 3 of \cite{ONP}, we obtain
\begin{equation}\label{ell1}
\ka =
\frac{\p(x)}{|\p(x)x-\aa|}e.
\end{equation}
 In this case all the caustic points coincide with the focus $\aa$ and $\ka$ depends on the point of reflection but not on a particular tangent direction at that point.
Note that
if the radial function is rescaled with a factor $\lambda >0$ and the eccentricity remains fixed then the second focus will be on the same axis but at
the distance $\lambda |\aa|$ from the first focus $\calO$.
The expression for $\ka$ remains invariant.

Similarly, for a one sheet of a two-sheeted hyperboloid of revolution with the revolution axis
passing through the foci, we obtain, by applying the corresponding expressions in section 3 of \cite{ONP},
\begin{equation}\label{ell2}\ka =
 -\frac{\p(x)}{|\p(x)x-\aa|}e.
\end{equation}

Just as in  the classical differential geometry,
we use the
quadratic forms $e$ and $\ka$ to define analogues of the principal curvatures and of the elementary symmetric functions of principal curvatures. For a hypersurface $R\in \M(\o)$ at a fixed point 
$x\in \o$ in an orthonormal basis such that $e_{ij}(x)=
\delta_{ij}$ the matrix $[\ka_{ij}]$ is symmetric and the roots $\lm_1,...,\lm_n$
of the polynomial equation
\begin{equation*}
P(\lm):=\det(a^i_j- \lambda \delta^i_j) =0,~\mbox{where} ~a^i_j:=e^{is}\ka_{sj},
\end{equation*}
are real. These roots will be called {\it principal} intensities.

For an integer $m,~1\leq m\leq n,$ define the {\it $m-$th intensity function} as the elementary
symmetric function
\begin{equation*}
S_m(\lambda) =
\sum_{1\leq i_1<...<i_m\leq n}\lambda_{i_1}\cdots \lm_{i_m}.
\end{equation*}
These functions are the coefficients of the polynomial
\[P(-\lm) = \lm^n - S_1\lm^{n-1}+...+(-1)^nS_n.\]
In particular,
\begin{equation}\label{intensity2}
S_n = \frac{\det [\ka_{ij}]}{\det [e_{ij}]}.
\end{equation}
It follows from (\ref{gigj}) and (\ref{density}) that
$|S_n|$ is the quotient of the volume forms defined by the form $\eh(\gR(x))$ and the metric
$e(x)$, that is, it is the quantity defined by (\ref{density}).
The analogue of the classical mean curvature is the {\it mean intensity} $(1/n)S_1$, where
\begin{equation}\label{mean}
S_1 = e^{ij}\ka_{ij}.
\end{equation}In view of (\ref{formb}), (\ref{H}) and definition of
$\ka$ we have
\[S_m(\lm_1(\p(x)),...,\lm_n\p((x)))=F_m(a^i_j(x))\equiv \]\[F_m(x, \p(x), \n_1 \p(x),...,\n_n \p(x), \n_{11}\p(x),...,\n_{nn}\p(x)), ~x \in \o,
\]
where $F_m$ is the sum of principal minors of $[a^i_j]$ of order $m$.

Fix some positive integer $m, ~1\leq m \leq m$. In analogy with the reflector problem we propose to study the  {\it $S_m$-reflector problem} for closed hypersurfaces in $\Rn$ concerned with determination of
 a closed
hypersurface $R\in \M(\Sn)$ such
that
\begin{equation}\label{Sm}
f(\gR(x))S_m(\lm_1(\p(x)),...,\lm_n(\p(x))) =  g(x), ~x \in \Sn,
\end{equation}
for  given positive functions $f$ and $g$ on $\Sn$; here $\p$ is the radial function of $R$. When $m=n$ this is the reflector problem described in the introduction.
For $m > 1$ these problems lead to fully nonlinear second order PDE's that have not yet been studied. The semilinear case when $n=1$ is treated in the next section.

Note that a positive solution $\p \in C^2(\Sn)$ of (\ref{Sm})
will always produce an embedded hypersurface in $\M(\Sn)$ with radial function $\p$.

It is worthwhile noting that if instead of the form $\ka$ the Schouten
tensor $\kappa(1 + |\n w|^2)/2$ is used (see Remark \ref{conformal}) and $\tilde{\lm}_1,...,\tilde{\lm}_n$ are its eigenvalues then the equation
\[S_m(\tilde{\lm}_1(w(x)),...,\tilde{\lm}_n(w(x))) =  c e^{-mw(x)} ~\mbox{on}~\Sn,\]
for some constant $c$, is the equation of the $S_m$-Yamabe problem  on $\Sn$ \cite{gursky:08}; here $w$ is as in Remark \ref{conformal}.

\section{Hypersurfaces with prescribed mean intensity}\label{ME_1}
It follows from (\ref{kappa0}) that $S_1$ in terms of $\p$ is given by
\begin{equation*}\label{me-0}
S_1(\lm_1(\p),...,\lm_n(\p)) =\frac{-\p \D \p +n\p^2+ 2|\n\p|^2 -(n/2)W_{\p}^2}{(1/2)W_{\p}^2}=:M[\p],
\end{equation*}
where $\D$ is the Laplace operator in the metric $e$.

Note that
for a sphere of radius $R$ with center at $\calO$ we have $S_1= n$ $\forall x \in \Sn$. For a hyperplane the mean intensity $S_1\equiv -n$. In this case the domain $\o$ is an open hemisphere. For a paraboloid of revolution $\ka \equiv 0$ and $S_1 \equiv 0$. The domain $\o$ in this case is $\Sn \setminus \{\xi\}$,
where $\xi$ is the axis of the paraboloid.
For an ellipsoid of revolution with one focus at
the origin $\calO$ and revolution axis passing through both foci $$S_1(\lm_1(\p(x)),...,\lm_n(\p(x)))=
\frac{n\p(x)}{|\p(x)x-\aa|}.$$

Setting $\p= 1/v$, we obtain a slightly simpler form of the operator $M$ above,
\begin{equation}\label{me-1}
M[1/v]= \frac{\D v +nv -nV}{V},~~\mbox{where}~~V= \frac{|\n v|^2+v^2}{2v}.
\end{equation}

The next proposition shows that there are no hypersurfaces in $\M(\Sn)$ with $S_1< n$ and $S_1> n ~\forall x \in \Sn$.
\begin{proposition}\label{estimate of S1}
Let $R\in \M(\Sn)$. Then there exist points on $\Sn$ where $S_1 \geq n$ and
$S_1 \leq n$.  Furthermore, the equality $S_1\equiv n$ is
attained only on concentric spheres centered at $\calO$.
\end{proposition}
\begin{proof} Suppose first that $S_1 > n~\forall x \in \Sn$. It follows from (\ref{me-1}) and the estimate $V \geq v/2$ that
\[
0=\int_{\Sn}\D v d\s =\int_{\Sn}(S_1+n)V d\s-n\int_{\Sn}v d\s \geq \frac{1}{2}\int_{\Sn}(S_1 -n)vd \s,
\]
where $d \s$ is the volume element on $\Sn$.
Thus, we arrived at a contradiction.

Suppose now that $S_1 < n~\forall x \in \Sn$. Let $x_0\in \Sn$ be a point where the $\min_{\Sn}v$ is attained. At $x_0$ we have: $\n v =0, ~V = v/2, ~\D v \geq 0$. Then
by (\ref{me-1}) at $x_0$ we have  $\D v = (S_1-n)(v/2)\geq 0$, which is impossible
if $S_1 < n$ on $\Sn$.

It remains to show that if $S_1\equiv n$ then $R$ is homothetic to $\Sn$.
To show this, note that in this case (\ref{me-1}) implies
\[
0=\int_{\Sn}\D v d\s =n \int_{\Sn}(2V-v) d\s.
\]
Since $2V \geq v $, we conclude that $2V = v$ and then $|\n v| = 0$. Hence,
$v = const$.
\end{proof}

Let $g:\Rn \rightarrow (0,\infty)$ be a given function. We write $\Rn \setminus\{\calO\}$ as $\Sn \times (0,\infty)$ and consider the  problem of finding a hypersurface $R\in \M(\Sn)$ defined by the radial function $\p: \Sn \rightarrow (0,\infty)$ and such that
\begin{equation}\label{me-2}
S_1(\lm_1(\p(x)),...,\lm_n(\p(x))) = \gba(x,\p(x)), ~x \in \Sn,
\end{equation}
where $\gba=ng$.
We have the following
\begin{theorem}\label{th1}
Let $g$ be a positive $C^{1,\a},~\a \in (0,1),$ function
in the annulus ${\calA }:=\{ x\in \Sn, ~\p \in [R_1,R_2]\}$,
where $0 < R_1 < R_2 < \infty$. Assume that $g$ satisfies the conditions:
\begin{equation}\label{th1-1}
(i) ~g(x,R_1) \geq 1~~\mbox{and}~~~(ii)~~g(x,R_2)\leq 1 ~~\forall x \in \Sn.
\end{equation}
Then there exists a hypersurface $R\in \M(\Sn)$ with radial function $\p \in C^{2,\a}(\Sn)$,
$\p(x) \in [R_1,R_2]~\forall x\in \Sn,$
satisfying the equation (\ref{me-2}).
\end{theorem}
\begin{proof}
Put $v=1/\p$. Then by (\ref{me-1}) we need to prove solvability of
the equation
\begin{equation}\label{me-7}
\D v +nv -nV =V\gba(x,1/v),~x \in \Sn.
\end{equation}
This is proved by applying the Leray-Schauder theorem on existence of fixed points to operator equations. In that we essentially follow the general scheme in O. Ladyzhenskaya and N. Uraltseva \cite{LU:68},  ch. IV, \S 10. However, the classes of functions we deal with here were not considered in \cite{LU:68} and we have to redo some of the steps and re-compute the degree of a certain map  that arises in our case.

Let
\begin{equation}\label{ineq0}
\Ca = \left \{w \in C^{1,\a}(\Sn)~|~\frac{1}{R_2}\leq w(x) \leq \frac{1}{R_1}~
\forall x \in \Sn\right \}
\end{equation}
and for
$v \in C^2(\Sn)$ and $w\in \Ca$ put
$$
\Dh v:=\D v+\frac{nv}{2},~~Q(x,w,\n w):= \frac{n|\n w|^2}{2w} + \frac{w^2 + |\n w|^2}{2w}\gba(x, 1/w),
$$
$$q(w):=\frac{nw^{1+\ep}\Rb^{\ep}}{2},~~Q^{\t}(x,w,\n w):=\t Q(x,w,\n w)+(1-\t)q(w),~\t\in [0,1],$$
where $\ep >0$ and $\Rb\in (R_1,R_2)$ are some fixed numbers.
Consider the family of problems
\begin{equation}\label{th3}
\Dh v = Q^{\t}(x,w, \n w),~ x\in \Sn,~ \t \in [0,1],~w \in \Ca.
\end{equation}
Note that when $\t=1$ and $w=v$ we obtain the equation (\ref{me-7}).

It is well known that the  two smallest eigenvalues of  -$\D$ on $\Sn$
are $0$ and $n$. Therefore, the uniformly elliptic operator $\Dh$ has a trivial kernel in $W^{1,2}(\Sn)$ and, consequently, in $C^{1,\a}(\Sn)$.
The function $g\in C^{1,\a}(\Sn \times [R_1, R_2])$. Hence, for any $w\in \Ca$ the right-hand side of (\ref{th3}) is in $C^{\a}(\Sn)$. Then, standard results on solvability of linear uniformly elliptic second order partial differential equations on $\Sn$ imply that the equation (\ref{th3}) has a unique solution $v_{\t}\in C^{2,\a}(\Sn)$ for any $w\in \Ca$ and $\t \in [0,1]$. Furthermore, by the Schauder estimate  for any solution of (\ref{th3})
\begin{equation}\label{LS1}
\parallel v_{\t}\parallel_{C^{2,\a}(\Sn)}\leq C \parallel \Dh v_{\t} \parallel_{C^{\a}(\Sn)},
\end{equation}
where $C>0$ is a constant depending on dimension $n$ and the coefficients of the standard metric of $\Sn$.
For domains in $\Rnn$ the Schauder inequality can be found, for example, in  \cite{LU:68}, ch. III, inequality (1.11). By applying it in  coordinate charts covering $\Sn$ (for example, the charts obtained with the stereographic projections from the North and South poles) and using a suitably large constant $C$, one obtains (\ref{LS1}). (The term $\max_{\Sn}|v_{\t}|$ usually included in the right hand side of (\ref{LS1}) is not needed here because  $\ker{\Dh}=\{0\}$.)
Thus, we have an operator
$T(w,\t):\Ca\times [0,1] \rightarrow C^{2,\a}(\Sn)$.

We want to apply the Leray-Schauder theorem  to the equation
\begin{equation}\label{LS0}
w=T(w,\t)
\end{equation}
 in the following setting. For a constant $A >0$ to be specified later, put
\begin{equation}\label{LSH}
U:=\left\{w \in \Ca \left|\right. \parallel w \parallel_{C^{1,\a}(\Sn)} < A+1\right \}
\end{equation}
and let $\Ub_1=\bar{U}\times[0,1]$.
Our goal is to show that $A$ can be chosen so that the following conditions hold:

(a) the set $U$ is connected and the map $T:\Ub_1\rightarrow C^{1,\a}(\Sn)$ is completely continuous,

(b) under the additional assumption that both inequalities in (\ref{ineq0}) are strict, the boundary of the set $U$ does not contain solutions of (\ref{LS0}) for all $\t\in [0,1]$, 

(c) for $\t=0$ the equation (\ref{LS0}) has a unique solution $w_0=1/\Rb$ and the linearized operator $(\Dh -Q^{\t}_w)$, where $Q^{\t}_w$ is the Fr\'{e}chet derivative evaluated on $w_0$ at $\t=0$ is invertible as a map from $C^{2,\a}(\Sn)$ to $C^{\a}(\Sn)$.

Under these circumstances, the Leray-Schauder degree of the operator $\mbox{Id}-T(\cdot,1)$ mapping $U$ to $0$ is defined and can be calculated to be $\pm 1$. Consequently, by the Leray-Schauder theorem there exists a $C^{2,\a}(\Sn)$ solution to (\ref{th3}) at $\t=1$. The additional assumption in (b) will be removed at the end of the proof.

Below, along with (\ref{LS0}), we will consider the equations
\begin{equation}\label{th3-1}
\Dh w= Q^{\t}(x,w \n w),~ \forall x\in \Sn,~\t\in [0,1].
\end{equation}
As it was already noted any $w\in \Ca$ substituted into the right
hand side of (\ref{th3-1}) for some $\t\in [0,1]$ gives a solution in $C^{2,\a}(\Sn)$. Thus, any $w\in C^2(\Sn)$ satisfying
(\ref{th3-1}) is, in fact, in $C^{2,\a}(\Sn)$. By construction, such $w$ 
also satisfies (\ref{LS0}) with the same $\t$. 

It is clear that the converse is also true. Namely, any $w\in \Ca$ satisfying (\ref{LS0}) for  some $\t\in [0,1]$ is in $\Ca\cap C^{2,\a}(\Sn)$ and  satisfies (\ref{th3-1}). Indeed, by construction, $T(w,\t)$ is a solution of (\ref{th3}) when such $w$ is inserted into the right hand side of (\ref{th3}). By Schauder's theorem $T(w,\t)\in C^{2,\a}(\Sn)$ and because of (\ref{LS0}) $w\in \Ca\cap C^{2,\a}(\Sn)$ and satisfies (\ref{th3-1}) (cf. \cite{LU:68}, p. 372).
This note is used below without further reminding.

Now, we prove (a). The connectedness of $U$ is clear as, in fact, $U$ is convex. To check that $T$ is completely continuous, we verify that 1) $T(w,\t)$ is continuous in $(w,\t)$ in $\Ub_1$ and continuous in $\t$ uniformly with respect to $w\in \Ub$ and 2)
for each fixed $\t \in [0,1]$ the map $T(w,\t)$ maps $\Ub$ into a compact set in $C^{1,\a}(\Sn)$.

To check 1) consider
$(w,\t), (w^{\prime},\t^{\prime}) \in \Ub_1$ and
the corresponding solutions $v$ and $v^{\prime}$. Put $\ov:=v^{\prime}-v$. Then
\begin{eqnarray}
\Dh \ov = \t^{\prime}[Q(x,w^{\prime}, \n w^{\prime})-Q(x, w, \n w)]\nonumber \\+(1-\t^{\prime})[q(w^{\prime})-q(w)]
+(\t^{\prime}-\t)[Q(x, w, \n w)]-q(w)].\label{LS2}
\end{eqnarray}
Using the interpolation $w^{s}:=sw^{\prime}+(1-s)w,~s \in [0,1]$, we obtain
\[
Q(x,w^{\prime}(x), \n w^{\prime}(x))-Q(x, w(x), \n w(x))=a^i(x)\n_i\ow(x) + a(x)\ow(x),\]
\[ q(w^{\prime}(x))-q(w(x))=b(x)\ow(x),
\]
where $\ow:=w^{\prime}-w,$ and
\[a_i(x)=\int_0^1\frac{\partial Q(x,w^s(x),\n w^s(x))}{\partial \n_i w^s}ds,~~a(x)= \int_0^1\frac{\partial Q(x,w^s(x),\n w^s(x))}{\partial w^s}ds,
\]
$$~~b(x)= \int_0^1\frac{\partial q(w^s(x))}{\partial w^s}ds.$$
It is clear that $a^i, a,b \in C^{\a}(\Sn),~i=1,...,n,$ and their $C^{\a}(\Sn)$ norms are bounded by a constant depending on $A$, the $\parallel \gba\parallel_{C^{1,\a}({\calA})}$ and $R_1,~R_2$.
Treating (\ref{LS2}) as a linear equation with respect to $\ov$ and applying (\ref{LS1}), we conclude that 1) is true.

As it was already noted $T(w,\t) \in C^{2,\a}(\Sn)$ for each $\t \in [0,1]$ and any $w  \in \Ub$. Since a set of functions bounded in the norm of $C^{2,\a}(\Sn)$ is compact
in $C^{1,\a}(\Sn)$, the operator $T(w,\t)$ maps $\Ub$
into a set compact in $C^{1,\a}(\Sn)$, that is, $T(w,\t)$ is a compact map from $\Ub_1$ into $C^{1,\a}(\Sn)$. This proves 2). Note that 1) and 2) hold with any $A < \infty$. This completes the proof of (a).

Next, we establish (b). We will need the following
\begin{lemma}\label{LSlemma} Suppose $w\in C^2(S^n)$ and  satisfies (\ref{th3-1}) for some $\t \in [0,1]$. Assume in addition that
\begin{equation}\label{LS6}
\frac{1}{R_2}\leq w(x)\leq \frac{1}{R_1}~\forall x \in \Sn.
\end{equation}
Then either $w\equiv 1/R_2$, or $w\equiv 1/R_1$, or
\begin{equation}\label{LS7}
\frac{1}{R_2}< w(x)< \frac{1}{R_1}~\forall x \in \Sn.
\end{equation}
\end{lemma}
The proof of this lemma will be given in the appendix.

{\it We now impose a temporary additional restriction on the         function $g$:}
\begin{equation}\label{LS-hyp}
g(x,R_1)\not \equiv 1  ~~\mbox{and}~~ g(x,R_2)\not \equiv 1~~\mbox{on}~~\Sn.
\end{equation}
Let us show that if (\ref{LS-hyp}) holds then
neither $1/R_1$ nor $1/R_2$ is
a solution of (\ref{th3-1}) 
for any $\t\in [0,1]$. Suppose $w(x)\equiv1/R_2$. Then for each $\t\in [0,1]$ and $\forall x \in \Sn$ we have
\begin{equation*}
D^{\t}(x):=\Dh \left (\frac{1}{R_2} \right )-Q^{\t}(x,\frac{1}{R_2}, 0) = \frac{n}{2R_2}\left [1 -\t g(x,R_2)-(1-\t)\left (\frac{\Rb}{R_2}\right)^{\ep}\right ].
\end{equation*}
Because of (\ref{LS-hyp}) and the inequality (ii) in (\ref{th1-1}) $g(\xb,R_2) < 1$ at some $\xb\in \Sn$. Since $\Rb/R_2 < 1$, we conclude that $D^{\t}(\xb) > 0$ for all $\t\in [0,1]$
and this proves our claim for $w\equiv 1/R_2$. The claim regarding $w\equiv1/R_1$ is proved similarly. 

Lemma \ref{LSlemma} implies now that under conditions (\ref{LS-hyp}) any $w\in C^2(S^n)$ satisfying  (\ref{th3-1}) for $\t\in [0,1]$ and 
\begin{equation}\label{LS9}
\frac{1}{R_2}\leq w(x)\leq \frac{1}{R_1}~\forall x \in \Sn
\end{equation}
is in fact such that
\begin{equation}\label{LS7-1}
\frac{1}{R_2}< w(x)< \frac{1}{R_1}~\forall x \in \Sn.
\end{equation}

Next, we check the applicability of the gradient estimates in \cite{LU:68},  ch. IV, Theorems 3.1 and 6.1, to solutions to (\ref{th3-1}) satisfying (\ref{LS9}).
For that we need to check
two conditions, the first of which (corresponding to (3.1) in ch. IV, \S 3) in our case reduces to positive
definiteness and uniform boundedness $\forall x \in \Sn$ of the quadratic form
$e^{ij}(x)\xi_i\xi_j,~\xi \in \Rnn$. This is  a consequence
of the properties of the metric $e$ on $\Sn$.
The  second condition (corresponding to (3.2) in ch. IV, \S 3) follows from the following estimate:  for all $x \in \Sn$ and any $w\in \Ca$ satisfying (\ref{LS9}) the inequality 
$$|\xi|(1+|\xi|)+|Q^{\t}(x,w,\xi)| \leq c_0(1+|\xi|)^2$$
holds for all $\tau \in [0,1]$ and all $\xi \in \Rnn$, where $|\xi|=(e^{ij}\xi_i\xi_j)^{1/2}$; the constant $c_0=c_0(n,R_1,R_2, \max_{\calA}g)$ and it is finite in a fixed coordinate atlas on $\Sn$.
It follows from Theorem 3.1 in \cite{LU:68},  ch. IV, that for any $\t\in [0,1]$ and any solution $w$
of (\ref{th3-1}) which is in $C^{2,\a}(\Sn)\cap \Ca$ the estimate
$$|\n w(x)| \leq c_1,$$
holds on $\Sn$; here, in each coordinate chart the constant $c_1$ depends on the same parameters
as $c_0$ above and on the distance from $x$ to the boundary of the chart.
By compactness of $\Sn$ we conclude that
$\parallel w \parallel_{C^1(\Sn)}$ is bounded by a constant depending only on
$n, R_1, R_2$ and $\parallel g \parallel_{C^1(\calA)}$. We keep the same notation, $c_1$, for that constant.

It remains to estimate the seminorm $|\n w|_{C^{\a}(\Sn)}$. We use for that a standard procedure; see, for example, \cite{LU:68},  ch. IV, \S 6, where this is done for domains in Euclidean space. Put $V:= \frac{w^2 + |\n w|^2}{2w}$,  fix an integer $1\leq s\leq n$ and differentiate covariantly the equation (\ref{th3-1}) with respect to the local variable $u^s$. Then, noting that $\n_sV = \frac{w^i\phi_{is}}{w}$, where $\phi_{is} = \n_{is} w + (w - V)e_{is}$, we get
\[
e^{ij}\n_s \n_{ij}w + \frac{n\n_s w}{2} =
\tau \left[\frac{n+\gba}{w} w^iq_{is} +
V(\n_s \gba+ \gba_w\n_sw ) - \frac{n\n_sw}{2}\right ]\]
\[+(1-\t)\frac{n(1+\ep)w^{\ep}\Rb^{\ep}\n_s w}{2};
\]
here, $\gba_w:=\del\gba/\del w$. By the Ricci identity
\[
 \n_s\n_{ij}w-\n_j\n_{is}w= e_{is}\n_jw -e_{ij}\n_sw.
\]
Putting $z := \n_s w$, we obtain
\[
\D z + \frac{2-n}{2}z -
\tau\left \{\frac{n+\gba}{w} w^i\n_iz +
\left [\frac{n+\gba}{w} (w - V )+V\gba_w-\frac{n}{2}\right]z
 \right \}\]
\[-(1-\t)\frac{n(1+\ep)w^{\ep}\Rb^{\ep}z}{2}= \t V \n_s \gba.
\]
(Note that $\n_i z = \partial_{si}w - \Gamma_{si}^k\n_kw$.)
This is a second order linear uniformly elliptic equation on $\Sn$ with respect to $z$.
For $w \in \Ca$ its coefficients are in $C^{\a}(\Sn)$ and by the Schauder estimate  (\cite{LU:68}, ch. III, inequality (1.11)) we have
$$\parallel z \parallel_{C^{2,\a}(\Sn)}\leq c_{2}(\parallel\t V \n_s \gba\parallel_{C^{\a}(\Sn)}+ \max_{\Sn}|z|),$$
where  $c_2=c_2(n,R_1,R_2,\max_{\calA}g, \parallel w \parallel_{C^1{\Sn}})$. Since $\parallel\t V \n_s \gba\parallel_{C^{\a}(\Sn)}$ and
$\max_{\Sn}|z|$ were already estimated through $n,R_1,R_2,\max_{\calA}g$, the last inequality implies that
$|z|_{C^{\a}(\Sn)}$ is bounded by a constant depending only on
$R_1, R_2, n$  and $\parallel g \parallel_{C^{1,\a}(\calA)}$.
Thus,  there exists a constant $c_3$, depending only on $n,R_1,R_2$, $\parallel g\parallel_{C^{1,\a}({\calA})}$ such that
the inequality
$$\parallel w\parallel_{C^{1,\a}(\Sn)} < c_3$$
holds for all
$w\in \Ca$ satisfying (\ref{LS0}) and each $\t \in [0,1]$.

Now, the constant $A$ in (\ref{LSH}) can be set equal to $c_3$. Then we conclude that under the restrictions (\ref{LS-hyp}) there are no solutions to (\ref{LS0}) on
the boundary of $U$.

Finally, we establish (c) and calculate the degree of $\mbox{Id} - T(\cdot,1)$. First we  consider  the equation (\ref{th3-1}) when $\t=0$:
\begin{equation}\label{LS4}
\Dh w=q(w) ~~\mbox{on}~~\Sn,~~w\in \Ub.
\end{equation}
A direct substitution shows that $w_0=1/\Rb$ is a solution of (\ref{LS4}). Let us show that this solution is unique.
Suppose $w^{\prime}\in \Ub$ is a solution of (\ref{LS4}) different from $1/\Rb$. At a point $x_{\max}\in \Sn$ where $w^{\prime}$ attains its maximum we have $\D w^{\prime}(x_{\max})\leq 0$ and then by (\ref{LS4}),
$$w^{\prime}(x_{\max})\leq \frac{1}{\Rb}.$$
Similarly, at a point $x_{\min}\in \Sn$ where $w^{\prime}$ attains its minimum $$w^{\prime}(x_{\min}) \geq \frac{1}{\Rb}.$$
Thus, the solution $w_0=1/\Rb$ of (\ref{LS4}) is unique.

Put $\Phi(w,\t):= \Dh w -Q^{\t}(x, w, \n w), ~w\in \Ca \cap C^{2,\a}(\Sn), ~\t \in [0,1].$ Then letting $w_s=w+sh,~h\in C^2(\Sn)$, calculating the derivative with respect to $s$ and setting $s=0$ we obtain
\[
\Phi_w(w,\t)(h):= \frac{d \Phi(w+sh, \t)}{ds}|_{s=0} =\Dh h - \t \frac{d Q(x, w+sh,
\n (w+sh))}{ds}|_{s=0}\]
\[-(1-\t)\frac{d q(w+sh)}{d s}|_{s=0}.
\]
The calculated weak derivative is in fact the Fr\'{e}chet derivative since it is uniformly continuous in $w$ in some $C^{2,\a}(\Sn)$ neighborhood
of $w$  and continuous in $h$ as a map from $C^{2,
\a}(\Sn)$ into $C^{
\a}(\Sn)$.

Evaluating the above expression on $w_0$ and $\t=0$, we get
\[\Phi_w(w_0,0)(h)=\D h-\frac{n\ep}{2}h.\]

Then
\[
\ker\Phi_w(w_0,0)= \{h\in C^2(\Sn)~|~\D h-\frac{n\ep}{2}h= 0 ~\mbox{on}~\Sn \}.
\]
Since we have chosen $\ep > 0$, it is easy to see that
$\ker\Phi_w(w_0,0)=\{0\}$. Standard results on linear elliptic partial differential equations imply that the map $\Phi_w(w_0,0):C^{2,\a} \rightarrow C^{\a}$ is an isomorphism.

By (a) and (b), the degree of the maps $\mbox{Id} - T(\cdot,t)$ into $0$ is defined for all $(w,t) \in \Ub_1$ and satisfying (\ref{LS7-1}). By standard results this degree is the same for all $t\in [0,1]$. Furthermore, since $\Phi_w(w_0,0):C^{2,\a}(\Sn) \rightarrow C^{\a}(\Sn)$ is an isomorphism, the derivative $\mbox{Id}-T_w(w_0,0)$ is invertible for all $w\in C^{\a}(\Sn)$ which are in a small $C^{\a}(\Sn)$
neighborhood of $w_0$ and satisfy (\ref{LS7-1}). On the other hand, $\mbox{deg} (\mbox{Id} - T(\cdot,0),U,0)=\mbox{deg} (\mbox{Id} - T_w(w_0,0),B,0)$, where $B$  is a ball in $C^{\a}(\Sn)$ with the center at $0$ and sufficiently small radius in the $C^{\a}(\Sn)$ norm (cf. \cite{nirenberg:01}, section 2.8). Since $\Phi_w(w_0,0)$ is invertible, its degree is $\pm 1$. Consequently, $\mbox{deg} (\mbox{Id} - T_w(w_0,0),B,0)=\pm 1$ and thus $\mbox{deg} (\mbox{Id} - T(\cdot,1),U,0)=\pm 1 \neq 0$. 
By the Leray-Schauder theorem, the equation  (\ref{th3}) has a fixed point
$w$ for $\tau=1$ in $\Ca$ and by the Schauder theorem
$w \in \Ca\cap C^{2,\a}(\Sn)$. This completes the proof of the theorem under  the restrictions (\ref{LS-hyp}).

If the restrictions (\ref{LS-hyp}) are not satisfied, that is,  $g(x,R_1)\equiv 1$ (or $g(x,R_2)\equiv 1$), then
a substitution of $v\equiv 1/R_1$ ($v\equiv 1/R_2$) into (\ref{me-7}) shows that $v\equiv 1/R_1$ ($v\equiv 1/R_2$) is a solution. Thus, the theorem is true also in these
cases.
\end{proof}

The next proposition deals with the question of uniqueness of  a solution found
in Theorem \ref{th1} and provides also
some additional information about such solutions.
\begin{proposition}
Suppose the conditions in Theorem \ref{th1} are satisfied and let $R^1$ and $R^2$ be two
hypersurfaces in $\M(\Sn)$ with radial functions $\p^1$ and $\p^2$ satisfying (\ref{me-2}).
If, in addition,
\begin{equation}\label{uniq}
\frac{\partial g}{\partial \p} \leq 0~\forall (x,\p) \in \Sn\times[R_1,R_2]
\end{equation}
then $\p^1(x) =C \p^2(x) ~\forall x\in \Sn$ and some constant $C > 0$.
Furthermore, each solution $\p\in C^2(\Sn)$ of (\ref{me-2}) such that
$R_1 \leq \p(x) \leq R_2 ~\forall x \in \Sn$ is either $\equiv R_1$ or
$\equiv R_2$ or
\begin{equation}\label{me-4}
R_1 < \p(x) < R_2 ~\forall x \in \Sn.
\end{equation}
\end{proposition}
\begin{proof}  Suppose $\p^2 > \p^1$ for some $x \in \Sn$. Put $\p^0:=C\p^2,$
where the constant $C\in (0,1)$ is chosen so that $\p^0(\ox)=\p^1(\ox)$ for some $\ox \in \Sn$
 and $\p^0(x)\leq\p^1(x)$ in some neighborhood $U\subset \Sn$  of $\ox$. Such $U$ is taken sufficiently small so that $g(x,\p^0(x))$ is defined for all $x\in U$.
Because $M$ is homogeneous of order zero in $\p$ and by (\ref{uniq}), we  have for all $x \in U$
\[M[\p^1] - \gba(x,\p^1(x))-M[\p^0]+\gba(x,\p^0(x)) \geq\]
\[ M[\p^1]-g(x,\p^1(x))-M[\p^2] +\gba(x,\p^2(x)) =0.\]
The operator $M[\p^t] -\gba(x,\p^t(x))$  is defined on $\p^t = (1-t)\p^0 + t\p^1,~\forall (x,t)\in U \times [0,1]$ and negatively uniformly elliptic. The left hand side of the  last inequality can be written as  
\[M[\p^1]-M[\p^0] - \gba(x,\p^1(x))+\gba(x,\p^0(x))=-M^{\prime}[\bar{\p}] - \bar{\p}(x)\int_0^1\frac{\partial \gba(x,\p^t(x))}{\partial \p}dt, \]
where $\bar{\p} = \p^1-\p^0$ and $M^{\prime}$ is a positively uniformly elliptic second order linear operator in $U$.
Thus,
$$M^{\prime}[\bar{\p}] + \bar{\p}(x)\int_0^1\frac{\partial \gba(x,\p^t(x))}{\partial \p}dt\leq 0~\forall x\in U.$$
Since $\bar{\p} \geq 0$ in $U$ and  $\bar{\p}(]\xb)=0$, the  strong maximum principle (see, \cite{Al_uniqueness-III}, Theorem B) implies
 $\p^0(x)-\p^1(x)\equiv 0$ in $U$ .
Consequently, the set $\{x \in \Sn~|~\p^0(x) = p^1(x)\}$ is open in $\Sn$. Since it is also closed, $\p^0(x)=\p^1(x) ~\forall x\in \Sn$ and therefore $\p^1=C\p^2$ with some $C \in (0,1)$. Reversing the roles of
$\p^2$ and $\p^1$, if necessary, we conclude that
$\p^1(x) =C \p^2(x) ~\forall x\in \Sn$ with some constant $C > 0$.

The last statement of the proposition follows from Lemma \ref{LSlemma}.
\end{proof}

\section{\protect\bf Appendix - Proof of Lemma \ref{LSlemma} }\label{appendix}
The claim in this lemma follows essentially from
Aleksandov's geometric (strong) form of the maximum principle.

Let $w\in C^2(\Sn)$ and satisfies (\ref{th3-1}) for some $\t \in [0,1]$. Assume also that (\ref{LS6}) holds. Suppose there exists some $x_0 \in \Sn$ such that
$w(x_0)=1/R_1$ and $w(x)\not \equiv 1/R_1$.
Consider
$$w^s(x):=(1-s)w(x)+\frac{s}{R_1},~(x,s)\in \Sn\times [0,1].$$
Observe that
$$\frac{1}{R_2}\leq w^s(x) \leq \frac{1}{R_1}~\forall (x,s)\in \Sn\times [0,1].$$
The operator $\Dh w^s - Q^{\t}(x,w^s,\n w^s)$ is defined and uniformly elliptic on $\Sn$ for all $s \in [0,1]$.

For $w^0=w$ and $w^1=1/R_1$ we have, taking into account that
$g(x,R_1) \geq 1$ and $\Rb > R_1$,
$$
\Dh w- Q^{\t}(x,w, \n w)-
\left [\Dh \left(\frac{1}{R_1}\right)-\t\frac{ng(x,R_1)}{2R_1}-
(1-\t)\frac{n \Rb^{\ep}}{2R_1^{1+\ep}}\right ]
$$
$$
 =-\frac{n}{2R_1}
\left [1-\t g(x,R_1)-(1-\t)\frac{\Rb^{\ep}}{R_1^{\ep}}\right ]\geq 0 ~\mbox{in}~\Sn.
$$
Since $w(x_0)=1/R_1$ and $w(x)\leq 1/R_1$ on $\Sn$ it follows from Theorem C in \cite{Al_uniqueness-III} that $w(x)\equiv 1/R_1$ in some neighborhood $\calK$ of $x_0$ on $\Sn$. By continuity of $w$ the equality $w(x)\equiv 1/R_1$ holds in $\bar{\calK}$. Thus, the set of points on $\Sn$ where $w(x)\equiv 1/R_1$ is open and closed on $\Sn$. Hence, $w(x)\equiv 1/R_1$ on $\Sn$.

Similarly, it is shown that $w\equiv 1/R_2$ if $w(x_1)=1/R_2$ at some $x_1\in \Sn$. This completes the proof of Lemma \ref{LSlemma}.

\bibliographystyle{plain}

\end{document}